\documentclass[11pt]{amsart}
\input epsf.sty
\usepackage{amscd,amssymb,verbatim}

\newcommand{\bd}{{\textrm{Bd}\,}}
\newcommand{\be}{\begin{enumerate}}

\newcommand{\bo}{\partial\,}

\newcommand{\chom}{{\mathcal Hom}}

\newcommand{\cn}{{\mathcal N}}

\newcommand{\com}{\complement}

\newcommand{\da}{\Delta}

\newcommand{\dz}{{\mathbb Z}}

\newcommand{\ee}{\end{enumerate}}

\newcommand{\ind}{{\text Ind}\,}

\newcommand{\lra}{\longrightarrow}

\newcommand{\nin}{\noindent}

\newcommand{\ovr}[1]{\overline{#1}}

\newcommand{\pr}{\noindent{\bf Proof. }}

\newcommand{\ra}{\rightarrow}

\newcommand{\sm}{\setminus}

\newcommand{\supp}{\text{\rm supp}\,}

\newcommand{\thom}{\text{\tt Hom}\,}
\newcommand{\thomp}{\text{\tt Hom}_+\,}
\newcommand{\ti}{\tilde }

\newcommand{\vt}{\vartheta}

\newcommand{\wti}{\widetilde }

\newcommand{\zz}{{{\mathbb Z}_2}}

\newtheorem{thm}{Theorem}[section]
\newtheorem{df}  [thm]{Definition}

\newtheorem{crl} [thm]{Corollary}
\newtheorem{prop}[thm]{Proposition}

\newtheorem{rem} [thm]{Remark}

\numberwithin{equation}{section}

\title {Topological obstructions to graph colorings}

\author{Eric Babson and Dmitry N. Kozlov}

\date{\noindent
\today\\[0.05cm]
            }

        \address{{\it First author:} Department of Mathematics,
          University of Washington, Seattle, U.S.A.}
        \email{babson@math.washington.edu}

        \address{{\it Second author, current:} Department of
          Mathematics, University of Bern, Switzerland}
        \email{kozlov@math-stat.unibe.ch}

        \address{{\it Second author, on leave from:} Department of
          Mathematics, Royal Institute of Technology, Stockholm,
          Sweden} \email{kozlov@math.kth.se}

\begin{document}

\begin{abstract}
  For any two graphs $G$ and $H$ Lov\'asz has defined a cell complex
  $\thom(G,H)$ having in mind the general program that the algebraic
  invariants of these complexes should provide obstructions to graph
  colorings. Here we announce the proof of a conjecture of Lov\'asz
  concerning these complexes with $G$ a cycle of odd length.

  \nin More specifically, we show that
\begin{quote}
  {\it If $\thom(C_{2r+1},G)$ is $k$-connected, then $\chi(G)\geq
    k+4$.}
\end{quote}

\nin Our actual statement is somewhat sharper, as we find obstructions
already in the non-vanishing of powers of certain Stiefel-Whitney
classes.

\end{abstract}

\maketitle


\section{Introduction}\label{sect_intr}

One of the central questions of graph theory is to find lower bounds
for the number of colors in an admissible coloring of the set of
vertices of a~graph (a coloring is called admissible if any 2
vertices, which are connected by an edge get different colors). The
minimal possible number of colors in an admissible coloring is called
the {\bf chromatic number} of the graph, and is denoted by $\chi(G)$.
In 1978 L.\ Lov\'asz solved the Kneser conjecture by finding geometric
obstructions of Borsuk-Ulam type to graph colorings.

\begin{thm}
  {\rm (Kneser-Lov\'asz, \cite{Knes,Lo}).}  {\it Let $\Gamma_{k,n}$ be
    the graph whose vertices are all $k$-subsets of $[n]$, and edges
    are all pairs of disjoint $k$-subsets; here $1\leq k\leq n/2$.
    Then $\chi(\Gamma_{k,n})=n-2k+2$.}
\end{thm}

The hard part was to show the inequality $\chi(\Gamma_{k,n})\geq
n-2k+2$.  Lov\'asz' idea to achieve that was to associate a~simplicial
complex $\cn(G)$, called the neighborhood complex, to an~arbitrary
graph $G$, and then use the connectivity information of the
topological space $\cn(G)$ to find obstructions to the colorability
of~$G$.

To define $\cn(G)$ we need some notations. For any graph $G$, let
$V(G)$ denote the set of its vertices, and $E(G)\subseteq V(G)\times
V(G)$ the set of its edges. Since all the graphs in this paper are
undirected, $(x,y)\in E(G)$ implies $(y,x)\in E(G)$. We allow our
graphs to contain loops. For a~graph $G$ we distinguish between looped
and unlooped complements, namely we let $\com G$ be the graph defined
by $V(\com G)=V(G)$, $E(\com G)=(V(G)\times V(G))\sm E(G)$, while
$\ovr{G}$ is the graph defined by $V(\ovr{G})=V(G)$,
$E(\ovr{G})=E(\com G)\sm\{(v,v)\,|\, v\in V(G)\}$. For a~graph $G$ and
$S\subseteq V(G)$ we denote by $G[S]$ the graph on the vertex set $S$
induced by $G$, that is $V(G[S])=S$, $E(G[S])=(S\times S)\cap E(G)$.

The complex $\cn(G)$ can then be defined by taking as its vertices all
non-isolated vertices of $G$, and as its simplices all the subsets of
$V(G)$ which have a~common neighbor; in other words, the maximal
simplices of $\cn(G)$ are $\cn(v)$, for $v\in V(G)$, where
$\cn(v)=\{w\in V(G)\,|\,(v,w)\in E(G)\}$ is the set of all neighbors
of~$v$.

\begin{thm} \label{lothm}
  {\rm (Lov\'asz, \cite{Lo}).}  {\it Let $G$ be a graph without loops, such that
    $\cn(G)$ is $k$-connected, then $\chi(G)\geq k+3$.}
\end{thm}

Numerous papers followed Lov\'asz' pioneering work, see \cite{Zi02}
for a good recent survey.  Generalizing the original line of thought,
for any pair of graphs $G$ and $H$, Lov\'asz has defined a~cell
complex $\thom(G,H)$.

Recall that a {\bf graph homomorphism} is a map $\phi:V(G)\rightarrow
V(H)$, such that $\phi\times\phi:E(G)\ra E(H)$, that is if $x,y\in
V(G)$ are connected by an edge, then $\phi(x)$ and $\phi(y)$ are also
connected by an edge. We denote the set of all homomorphisms from $G$
to $H$ by $\chom(G,H)$.

\begin{df} \label{dfhom}
  $\thom(G,H)$ is a polyhedral complex whose cells are indexed by all
  functions $\eta:V(G)\rightarrow 2^{V(H)}\setminus\{\emptyset\}$,
  such that if $(x,y)\in E(G)$, for any $\tilde x\in\eta(x)$ and
  $\tilde y\in\eta(y)$ we have $(\tilde x,\tilde y)\in E(H)$.

  The closure of a~cell $\eta$ consists of all cells indexed by
  $\ti\eta:V(G)\rightarrow 2^{V(H)}\setminus\{\emptyset\}$, which
  satisfy $\ti\eta(v)\subseteq\eta(v)$, for all $v\in V(G)$.
\end{df}

In particular, the set of vertices of $\thom(G,H)$ is precisely
$\chom(G,H)$.

\vspace{5pt}

As the Proposition~\ref{propn} shows, already the complexes
$\thom(K_2,G)$ contain the homotopy type information of the
neighbourhood complexes. This observation makes one try to look for
obstructions coming from maps from other graphs than $K_2$.

In this paper we announce the proof of the following conjecture of
Lov\'asz, which is a natural next step, after the Theorem~\ref{lothm},
in the general program of finding topological obstructions to the
existence of the graph homomorphisms.

\begin{thm} \label{lovcon}
  Let $G$ be a graph, and let $r,k\in\dz$, such that $r\geq 1$, $k
  \geq -1$.  If $\thom(C_{2r+1},G)$ is $k$-connected, then
  $\chi(G)\geq k+4$.
\end{thm}
Here, $C_{2r+1}$ is a~$(2r+1)$-cycle, i.e.\ a~graph defined by
$V(C_{2r+1})=\dz/(2r+1)\dz$,
$E(C_{2r+1})=\{([x],[x+1]),([x+1],[x])\,|\,x\in\dz/(2r+1)\dz\}$.

To state our main result we need more notations. Let $X$ be a~CW
complex with a free $\zz$-action. By the general theory of principal
bundles, there exists a~$\zz$-equivariant map $\tilde w:X\ra
S^{\infty}$, where $\zz$ acts on $S^{\infty}$ by the antipodal map,
and the induced map $w:X/\zz\ra{\mathbb R\mathbb P}^\infty$ is unique
up to homotopy. This in turn induces a~canonical algebra homomorphism
$w^*:H^*({\mathbb R\mathbb P}^\infty;\zz)\ra H^*(X/\zz;\zz)$. Let
$z\in H^1({\mathbb R\mathbb P}^\infty;\zz)$ denote the nontrivial
cohomology class, then $H^*({\mathbb R\mathbb
  P}^\infty;\zz)\simeq\zz[z]$ as graded algebras. We denote $w^*(z)\in
H^1(X/\zz;\zz)$ by $\varpi_1(X)$, and call it the {\bf first
  Stiefel-Whitney class} of the $\zz$-space~X. Clearly,
$\varpi_1^k(X)=w^*(z^k)$.

Let $\zz$ act on $C_{2r+1}$ by mapping $[x]\in\dz/(2r+1)\dz$ to
$[-x]$, and let $\gamma\in\chom(C_{2r+1},C_{2r+1})$ denote the
corresponding graph homomorphism. Furthermore, let $\zz$ act on
$K_m$ for $m\geq 2$, by swapping the vertices 1 and 2 and fixing
the vertices $3,\dots,m$; here, $K_m$ is the graph defined by
$V(K_m)=[m]$, $E(K_m)=\{(x,y)\,|\,x,y\in [m], x\neq y\}$. These
induce free $\zz$-actions on $\thom(C_{2r+1},G)$ and
$\thom(K_m,G)$.

\begin{thm} \label{thmmain}
 Let $G$ be a graph.
\begin{enumerate}
\item[(a)] If $\varpi_1^k(\thom(C_{2r+1},G))\neq 0$, then $\chi(G)\geq k+3$.
\item[(b)] If $\varpi_1^k(\thom(K_m,G))\neq 0$, then $\chi(G)\geq k+m$.
\end{enumerate}
\end{thm}

Since we know that for a $k$-connected space $X$ with a~free
$\zz$-action $\varpi_1^{k+1}(X)\neq 0$, the Theorem \ref{thmmain}
implies the Theorem \ref{lovcon} as well as the following corollary.

\begin{crl}
  Let $G$ be a graph, and let $m,k\in\dz$, such that $m\geq 2$, $k\geq
  -1$. If $\thom(K_m,G)$ is $k$-connected, then $\chi(G)\geq k+m+1$.
\end{crl}

\nin {\bf Acknowledgments.} We would like to thank L\'aszl\'o
Lov\'asz and P\'eter Csorba for insightful discussions. The second
author acknowledges support by the University of Washington,
Seattle, the Swiss National Science Foundation, and the University
of Bern.

\section{Homotopy type of the  $\thom$ complexes}

First, we remark two properties of the $\thom$ complexes:
\begin{enumerate}
\item[(1)] Cells of $\thom(G,H)$ are direct products of simplices.
  More specifically, each $\eta$ as in the Definition~\ref{dfhom} is
  a~product of $|V(G)|$ simplices, having dimensions $|\eta(x)|-1$,
  for $x\in V(G)$.
\item[(2)] $\thom(H,-)$ is a covariant, while $\thom(-,H)$ is a
  contravariant functor from {\bf Graphs} to {\bf Top}, where {\bf
    Graphs} is a category having graphs as objects, and graph
  homomorphisms as morphisms. If $\phi\in\chom(G,G')$, then we shall
  denote the induced topological maps as
  $\phi^H:\thom(H,G)\ra\thom(H,G')$ and
  $\phi_H:\thom(G',H)\ra\thom(G,H)$.
\end{enumerate}

\begin{prop} \label{propn}
$\thom(K_2,G)$ is homotopy equivalent to $\cn(G)$.
\end{prop}

If the Proposition \ref{propn} coupled with the Theorem~\ref{lothm}
were to be interpreted as that the obstructions to colorability found
by Lov\'asz in 1978 stem from the graph $K_2$, then the idea behind
the Lov\'asz Conjecture could be thought of as that the next natural
class of obstructions should come from the odd cycles, $C_{2r+1}$.

The next proposition is the crucial step in the proof of the
Theorem \ref{thmmain}(b).

\begin{prop} \label{pr_chom}
$\thom(K_m,K_n)$ is homotopy equivalent to a wedge of $(n-m)$-dimensional
spheres.
\end{prop}

Next, we need a technical Quillen-type lemma. For any small category
$C$ (in particular a finite poset) we denote by $\Delta(C)$ the
realization of the nerve of that category. For any finite poset $P$,
we let $P^{op}$ denote the finite poset which has the same set of
elements as $P$, but the opposite partial order. Also, for any finite
poset $P$, whenever the subset of the elements of $P$ is considered as
a poset, the partial order is taken to be induced from $P$.

\begin{prop} \label{prop_ABC}
Let $\phi:P\ra Q$ be a map of finite posets. Consider a list of
possible conditions on $\phi$.

\nin {\rm Condition $(A)$.} For every $q\in Q$, $\da(\phi^{-1}(q))$ is
contractible.

\nin {\rm Condition $(B)$.} For every $p\in P$ and $q\in Q$ with
$\phi(p)\geq q$ the poset $\phi^{-1}(q)\cap P_{\leq p}$ has a
maximal element.

 \nin {\rm Condition $(B^{op})$.} Let
$\phi^{op}:P^{op}\ra Q^{op}$ be the poset map induced by $\phi$.
We require that $\phi^{op}$ satisfies Condition~$B$.

Then
\begin{enumerate}
\item [(1)] If $\phi$ satisfies $(A)$ and either $(B)$ or $(B^{op})$,
  then $\phi$ is a homotopy equivalence.

\item [(2)] If $\phi$ satisfies $(B)$ and $(B^{op})$, and $Q$ is
  connected, then for any $q,q'\in Q$ we have
  $\da(\phi^{-1}(q))\simeq\da(\phi^{-1}(q'))$.  Furthermore, we have
  a~fibration homotopy long exact sequence:

\begin{equation}
\dots\lra\pi_i(\da(\phi^{-1}(q)))\lra\pi_i(\da(P))\lra\pi_i(\da(Q))\lra\dots
\end{equation}
\end{enumerate}
\end{prop}
\pr The argument is based on using the Quillen's Theorems A and B, see
\cite[pp.\ 85,89]{Qu}, for the induced map $\bd\phi:\bd P\ra\bd Q$,
here $\bd$ denotes the barycentric subdivision, that is the poset of
all the chains in the given poset. The details will appear elsewhere.
\qed

\begin{prop}
  If $G$ and $H$ are graphs and $u$ and $v$ are vertices of $G$, such
  that $N(v)\subseteq N(u)$, then $i:G-v\hookrightarrow G$ induces a
  homotopy equivalence $i_H:\thom(G,H)\ra\thom(G-v,H)$.
\end{prop}
\pr Apply the Proposition \ref{prop_ABC} (1) for the cellular map
$i_H:\thom(G,H)\ra\thom(G-v,H)$, which forgets the colors of $v$.
\qed

\begin{crl} \label{homtree}
  If $T$ is a~tree with at least one edge, then the map
  $i_{K_n}:\thom(T,K_n)\ra\thom(K_2,K_n)$ induced by any inclusion
  $i:K_2\hookrightarrow T$ is a~homotopy equivalence. Furthermore, if
  $F$ is a~forest, then $\thom(\ovr{F},K_n)\simeq\thom(K_m,K_n)$, where
  $m$ is the maximal cardinality of an~independent set in~$F$.
\end{crl}

\section{$\thomp$ and filtrations}

For a finite graph $H$, let $H_+$ be the graph obtained from $H$ by
adding an extra vertex, called the base vertex, and connecting it by
edges to all the vertices of $H_+$ including itself.

\begin{df}
  Let $G$ and $H$ be two graphs. The simplicial complex $\thomp(G,H)$
  is defined to be the link in $\thom(G,H_+)$ of the homomorphism
  mapping every vertex of $G$ to the base vertex in $H_+$.
\end{df}

So $\thomp(G,H)$ is just like $\thom(G,H)$ with the difference that we
also allow empty lists of colors, i.e., the cells are functions
$\eta:V(G)\rightarrow 2^{V(H)}$, which also makes it simplicial. We
remark that $\thomp(H,-)$ is a~covariant functor from {\bf Graphs} to
{\bf Top}.

\begin{prop}
$\thomp(G,H)$ is isomorphic to the independence complex of
$G\times\com H$. In particular, $\thomp(G,K_n)$ is isomorphic to
$\ind(G)^{*n}$, where $*$ denotes the simplicial join.
\end{prop}

The chain complex of $\thomp(G,K_n)$ has a natural filtration. For
each simplex of $\thomp(G,K_n)$, $\eta:V(G)\ra 2^{V(K_n)}$, define the
support of $\eta$ to be $\supp\eta=V(G)\setminus\eta^{-1}(\emptyset)$.
In fact, the same definition for $\supp$ holds for any $\thomp(G,H)$.
One concise way to phrase it is to simply consider the map
$t^G:\thomp(G,H)\ra\thomp(G,\com K_1)\simeq\Delta_{|V(G)|-1}$ induced by
the homomorphism $t:H\ra\com K_1$.

We can now filter $C_*(\thomp(G,K_n))$ by the cardinalities of the
support. Namely,
\[
F_{-1}=0,\,F_0=\langle|\supp\eta|=1\rangle,\dots,
F_k=\langle|\supp\eta|\leq k+1\rangle,\dots,
\]
that is $F_k$ is the chain subcomplex of $C_*(\thomp(G,K_n))$
generated by all simplices whose support set has cardinality at most
$k+1$.

Clearly, for any $k$,
\[
F_k/F_{k-1}=C_*(\coprod_{S\subseteq V(G), |S|=k+1}
\thom(G[S],K_n))[k],
\]
where the brackets $[-]$ denote shifting. In particular,
$F_{|V(G)|-1}=C_*(\thomp(G,K_n))$, and
$F_{|V(G)|-1}/F_{|V(G)|-2}=C_*(\thom(G,K_n))[|V(G)|-1]$.

\vspace{5pt}

Next, we describe a natural filtration on
$C_*(\thomp(C_{2r+1},K_n)/\zz;\dz)$. Let
$c,a_1,\dots,a_r,b_1,\dots,b_r$ denote the vertices of $C_{2r+1}$ so
that $\gamma(a_i)=b_i$, for any $i\in[r]$, and
$(c,a_1),(a_i,a_{i+1})\in E(C_{2r+1})$, for any $i\in[r-1]$. Identify
$V(C_{2r+1})$ with the vertices of an abstract simplex $\Delta_{2r}$
of dimension $2r$.

We subdivide $\Delta_{2r}$ by adding $r$ more vertices, denoted
$c_1,c_2,\dots,c_r$, and defining a~new abstract simplicial complex
$\tilde\Delta_{2r}$ on the set
$\{c,a_1,\dots,a_r,b_1,\dots,b_r,c_1,\dots,c_r\}=V(\tilde\Delta_{2r})$.
The simplices of $\tilde\Delta_{2r}$ are all the subsets of
$V(\tilde\Delta_{2r})$ which do not contain the subset $\{a_i,b_i\}$,
for any $i\in[r]$.

One can think of this new complex $\tilde\Delta_{2r}$ as the one
obtained from $\Delta_{2r}$ by representing it as a~join
$\{c\}*[a_1,b_1]*\dots*[a_r,b_r]$, inserting one additional vertex
into the middle of each $[a_i,b_i]$, and then taking the join of
$\{c\}$ and the subdivided intervals.  For
$\tilde\sigma\in\tilde\Delta_{2r}$ we define
$\vartheta(\tilde\sigma)\in\Delta_{2r}$ by $\vartheta(\tilde\sigma)=
(\tilde\sigma\sm\{c_1,\dots,c_r\})\cup
\bigcup_{c_i\in\ti\sigma}\{a_i,b_i\}$.

$\tilde\Delta_{2r}$ has an additional property: if a~simplex of
$\tilde\Delta_{2r}$ is $\gamma$-invariant, then it is fixed pointwise.
This allows us to introduce a~simplicial structure on
$\tilde\Delta_{2r}/\zz$ by taking the orbits of the simplices of
$\tilde\Delta_{2r}$ as the simplices of $\tilde\Delta_{2r}/\zz$.

Let us now describe a chain complex $\wti C_*(\thomp(C_{2r+1},K_n))$,
which comes from a different triangulation of the topological space
$\thomp(C_{2r+1},K_n)$. The chain complex consists of vector spaces
over~$\zz$. The generators are pairs $(\eta,\sigma)$, where
$\eta\in\thomp(C_{2r+1},K_n)$, and $\sigma\in\tilde\Delta_{2r}$, such
that $\vartheta(\sigma)=\supp\eta$. Such a~pair corresponds to the
cell $\eta\cap\supp^{-1}(\sigma)$. The codimension~1 boundary of
$(\eta,\sigma)$ is the sum of the following generators:
\begin{enumerate}
\item[(1)] $(\ti\eta,\sigma)$, if $\supp\ti\eta=\supp\eta$,
and $\ti\eta\in\bo\eta$;
\item[(2)] $(\eta,\ti\sigma)$, if $\ti\sigma\in\bo\sigma$,
and $\vt(\sigma)=\vt(\ti\sigma)$;
\item[(3)] $(\ti\eta,\sigma\sm\{x\})$, if $x\in V(\ti\Delta_{2r})$,
$Y=\vt(\{x\})\sm\vt(\sigma\sm\{x\})\neq\emptyset$, where $\ti\eta$ is
obtained from $\eta$ by setting $\eta$ to be $\emptyset$ on the vertex
set $\eta^{-1}(Y)$, and where we require that all the values of
$\eta|_{\eta^{-1}(Y)}$ have cardinality~1.
\end{enumerate}
Here $\bo$ denotes the codimension~1 boundary. The degree of
$(\eta,\sigma)$ in $\wti C_*(\thomp(C_{2r+1},K_n))$ is given by
$$\deg(\eta,\sigma)=|\sigma|-1+\sum_{v\in V(C_{2r+1}),\,|\eta(v)|\neq
  0}(|\eta(v)|-1)=\deg\eta+|\sigma|-|\vt(\sigma)|.$$

$\zz$ acts on $\wti C_*(\thomp(C_{2r+1},K_n))$ and we let $\wti
C_*^\zz(\thomp(C_{2r+1},K_n))$ denote its subcomplex consisting of the
invariant chains. By construction of the subdivision, $\wti
C_*^\zz(\thomp(C_{2r+1},K_n))$ is a chain complex for a~triangulation
of the space $\thomp(C_{2r+1},K_n)/\zz$.

Just like before, we consider the natural filtration $(\wti F_{-1}
\subseteq\wti F_{0}\subseteq\dots)$ on the chain complex $\wti
C_*^\zz(\thomp(C_{2r+1},K_n))$ by the cardinality of~$\sigma$. This
time, for any $k$
\[
\wti F_k/\wti F_{k-1}=C_*(\coprod_\sigma(\thom(G[\vt(\sigma)],K_n)/\zz)\,\,
\coprod_\tau\thom(G[\vt(\tau)],K_n))[k],
\]
where the first coproduct is taken over all
$\sigma\in\tilde\Delta_{2r}$ which are $\zz$-invariant, and the
second coproduct is taken over all $\tau\in\tilde\Delta_{2r}$
which are not $\zz$-invariant. In addition $|\sigma|=|\tau|=k+1$
is required.

\section{The spectral sequences and the proof of the Lov\'asz Conjecture}

We shall show that
\begin{equation} \label{eqchar}
(\ti\iota_{K_n})_{n-2}:H_{n-2}(\thom(C_{2r+1},K_n)/\zz)\ra
H_{n-2}(\thom(K_2,K_n)/\zz),
\end{equation}
is a 0-map, where the homology is taken with coefficients in $\zz$.
Here $\iota:K_2\hookrightarrow C_{2r+1}$ is either of the two
$\zz$-equivariant inclusion maps which take the vertices of $K_2$ to
$\{a_r,b_r\}$. It induces a $\zz$-equivariant map
$\iota_{K_n}:\thom(C_{2r+1},K_n)\ra\thom(K_2,K_n)$, and hence also the
quotient map $\ti\iota_{K_n}:\thom(C_{2r+1},K_n)/\zz\ra
\thom(K_2,K_n)/\zz$.

Note that if $(\ti\iota_{K_n})_{n-2}$ is a 0-map, then so is
$$(\ti\iota_{K_n})^{n-2}:H^{n-2}(\thom(K_2,K_n)/\zz;\zz)\ra
H^{n-2}(\thom(C_{2r+1},K_n)/\zz;\zz),$$ since we are working over
a~field. By the functoriality of the Stiefel-Whitney classes, this
implies that $\varpi_1^{n-2}(\thom(C_{2r+1},K_n))=0$, which is
a~crucial fact for our proof of Lov\'asz Conjecture. We shall
sketch the full computation for the case $n$ is odd. For even $n$
the Lov\'asz Conjecture follows already from \eqref{stern}, see
the discussion following~\eqref{stern}.

The $E^1$-tableau of the first spectral sequence is given by
$E_{d,s}^1=H_d(F_s,F_{s-1})$. From our remarks above concerning this
filtration it follows that each $E_{d,s}^1$ is a direct sum of the
appropriate homology groups of $\thom(H,K_n)$, where $H$ is a subgraph
of $C_{2r+1}$. Since all proper subgraphs of $C_{2r+1}$ are forests,
we know the first tableau, except for the $(2r)$th row, by the
Corollary~\ref{homtree}. It can be shown (we omit the technical
details) that the only part of the spectral sequence which is relevant
for the computation of the homology groups of $\thom(C_{2r+1},K_n)$ up
to dimension $n-2$, consists of two chain complexes:
\[
D_*^0 : E_{0,0}^1\ra E_{1,1}^1\ra E_{2,2}^1\ra\dots\ra E_{2r-1,2r-1}^1
\]
and
\[
D_*^1 : E_{n-2,0}^1\ra E_{n-1,1}^1\ra E_{n,2}^1\ra\dots\ra
E_{n-2+2r-1,2r-1}^1,
\]
where all the differentials are $d^1$ from the spectral sequence. To
start with, $D_*^0$ is isomorphic to the chain complex of the
$(2r-1)$-skeleton of the $(2r)$-simplex $\Delta_{2r}$, thus having
only the homology in dimensions 0 and $2r-1$. We can conclude that
$E_{0,0}^2=\dz$, while $E_{1,1}^2=E_{2,2}^2=\dots=E_{2r,2r}^2=0$,
hence $E_{2r,2r}^1=\dz$.

The analysis of $D_*^1$ is a little more interesting. By taking the
induced subgraphs, we identify subsets of $V(C_{2r+1})$ with
collections of arcs on a~circle (some arcs of length 1), such that
every two arcs are separated by at least one vertex. The detailed
analysis using the Corollary~\ref{homtree} shows that the generators
of the vector spaces in $D_*^1$ can be indexed with pairs $(S,A)$,
where $S\subseteq V(C_{2r+1})$, and $A$ a~connected component of
$C_{2r+1}[S]$ with at least one edge. We call such components {\it
  arcs}, and we call this specific $A$ a~{\it marked arc}. The
dimension of $(S,A)$ is $|S|$. The boundary of $(S,A)$ is the sum
(with appropriate signs) over all pairs obtained by deleting an
element of $S$. One can work out that, if the deleted element is not
in $A$, then the arc remains the same, while a~deleted element which
is in $A$ contributes to the boundary a sum of two pairs, depending on
which piece of $A$ is taken as a new arc, where the pieces of length 1
give 0 contribution.

We can now filter $D_*^1$ by the length of the marked arc $A$. By the
previous argument, the relative chain complexes will split into
subcomplexes each corresponding to a particular marked arc. All these
subcomplexes are acyclic (they are reduced chain complexes of
simplices) except for two, namely, those corresponding to
$S=V(C_{2r+1})\sm\{c\}$ and $S=V(C_{2r+1})\sm\{a_r,b_r\}$.

The complete computation of the first differential shows that
\begin{equation} \label{stern}
H_i(\thom(C_{2r+1},K_n);\zz)=\begin{cases}
  \zz,&\text{ if }i=0,n-3,n-2,\\
  0,&\text{ if }1\leq i\leq n-4.
\end{cases}
\end{equation}

Additionally, one can see that $H_i(\thom(C_{2r+1},K_n);\dz)=0$,
for $1\leq i\leq n-4$, and
$H_{n-2}(\thom(C_{2r+1},K_n);\dz)=H_{n-3}(\thom(C_{2r+1},K_n);\dz)=\dz$,
for $n$ odd, while $H_{n-2}(\thom(C_{2r+1},K_n);\dz)=0$,
$H_{n-3}(\thom(C_{2r+1},K_n);\dz)=\zz$, for $n$ even.

We remark that from this computation Lov\'asz Conjecture follows
easily in the case $n$ is even. The essence is to show that there
is no $\zz$-map of $S^{n-2}$ into $\thom(C_{2r+1},K_n)$. If this
map existed, then, combined with the canonical map
$\thom(C_{2r+1},K_n)\ra\thom(K_2,K_n)$, induced by
$K_2\hookrightarrow C_{2r+1}$, it would yield a $\zz$-map of
degree 0 of $S^{n-2}$ into itself, which is impossible.

We assume for the rest of the paper that $n$ is odd. In this case
one can show that the $\zz$-action on
$H_{n-3}(\thom(C_{2r+1},K_n);\dz)$ is always multiplication
by~$-1$.

The $E^1$-tableau of the second spectral sequence is given by
$E_{d,s}^1=H_d(\wti F_s,\wti F_{s-1})$. This time each $E_{d,s}^1$
splits as a~vector space over $\zz$ into direct sums of
$H_i(\thom(H,K_n);\zz)$, where $H$ is an~induced subgraph of
$C_{2r+1}$, and of $\thom(H,K_n)/\zz$, where $H$ is a~$\zz$-invariant
induced subgraph of $C_{2r+1}$. The detailed analysis of the first
tableau of this spectral sequence is rather technical and will appear
elsewhere. We sketch here the basic steps of the argument.

\nin (1) We single out the generator $\rho$ in $E_{n-2+r,r}^1$, which
is the contribution of $\{c,c_1,\dots,c_r\}$. Let $\wti E_{n-2+r,r}^1$
denote the part of $E_{n-2+r,r}^1$ which is spanned by the rest of the
generators.

\nin (2) Analysis of the boundary of $\rho$ shows that if
$(\ti\iota_{K_n})_{n-2}$ is not a~0-map, then rank of
$d^1:E_{n-2+r,r}^1\ra E_{n-3+r,r-1}^1$ is one higher than the rank of
$\wti d^1:\wti E_{n-2+r,r}^1\ra E_{n-3+r,r-1}^1$. This is because the
boundary of other generators of $E_{n-2+r,r}^1$ does not contain the
generator in $E_{n-3+r,r-1}^1$ indexed by $\{c_1,\dots,c_r\}$. Indeed,
the only place for it to be would be in the boundary of the generators
indexed by $\{x,c_1,\dots,c_r\}$, where $x=a_i$ or $x=b_i$, but
omitting $x$ is a double covering map, which is a 0-map on the
$\zz$-homology.

\nin(3) The differential $d^1:E_{n-2+r,r+1}^1\ra E_{n-3+r,r}^1$ is
a~0-map, because $q_{n-3}$ in \eqref{eqquot} below is a~0-map.

\nin(4) The chain complex
\[
E_{n-2+r,r+1}^1\ra E_{n-1+r,r+2}^1\ra\dots\ra E_{n-3+2r,2r}^1
\]
is a chain complex of ${\mathbb R}{\mathbb P}^{r-1}$ over $\zz$. Hence
$E_{n-2+r,r+1}^2=\zz$, and so $d^2:E_{n-2+r,r+1}^2\ra E_{n-3+r,r-1}^2$
has rank~1.

\nin(5) By means of technical computations we can show that the sequence
\[
E_{n-4+r,r-2}^1\ra E_{n-3+r,r-1}^1\ra \wti E_{n-2+r,r}^1
\]
has homology $\zz$ in the middle term.  Altogether, this shows the
equation \eqref{eqchar}.

It remains to see that
\begin{equation}\label{eqquot}
q_{n-3}:H_{n-3}(\thom(C_{2r+1},K_n);\zz)\ra
H_{n-3}(\thom(C_{2r+1},K_n)/\zz;\zz).
\end{equation}
is a 0-map.

\noindent
{\bf Proof of \eqref{eqquot}.} First let us see \eqref{eqquot} over
integers.  Take $\zeta\in H_{n-3}(\thom(C_{2r+1},K_n);\dz)$. By our
previous computations $\gamma^{K_n}(\zeta)=-\zeta$, this means that
$q_{n-3}(\gamma^{K_n}(\zeta))= -q_{n-3}(\zeta)$. On the other hand,
$q_{n-3}$ commutes with the $\zz$-action, that is
$q_{n-3}(\gamma^{K_n}(\zeta))=q_{n-3}(\zeta)$, which yields
$q_{n-3}(\zeta)=0$.

Second, by the universal coefficient theorem the map
$$\tau:H_{n-3}(\thom(C_{2r+1},K_n);\dz)\otimes\zz\ra
H_{n-3}(C_*(\thom(C_{2r+1},K_n);\dz)\otimes\zz)$$
is injective and functorial. In our concrete situation, this map is
also surjective, hence the claim follows from the following diagram:
$$\begin{CD}
H_{n-3}(\thom(C_{2r+1},K_n);\dz)\otimes\zz @>\text{0-map}>>
H_{n-3}(\thom(C_{2r+1},K_n)/\zz;\dz)\otimes\zz \\
        @V\tau V\text{iso}V           @VVV\\
H_{n-3}(\thom(C_{2r+1},K_n);\zz)   @>q_{n-3}>>
H_{n-3}(\thom(C_{2r+1},K_n)/\zz;\zz) \qed \\
\end{CD}
$$

\vskip7pt

\noindent
{\bf Proof of the Lov\'asz Conjecture.}

\vskip5pt

\nin If the graph $G$ is $(k+3)$-colorable, then there exists
a~homomorphism $\phi:G\ra K_{k+3}$. It induces a $\zz$-equivariant
map $\phi^{C_{2r+1}}:\thom(C_{2r+1},G)\ra
\thom(C_{2r+1},K_{k+3})$. Since the Stiefel-Whitney classes are
functorial and $\varpi_1^{k+1}(\thom(C_{2r+1},K_{k+3}))=0$, the
existence of the map $\phi^{C_{2r+1}}$ implies that
$\varpi_1^{k+1}(\thom(C_{2r+1},G))=0$, which is a contradiction to
the assumption of the theorem. \qed


\begin{thebibliography}{00}

\bibitem{Knes} M.\ Kneser, Aufgabe 300, {\it Jber.\ Deutsch.\ Math.-Verein.\ }
{\bf 58} (1955).

\bibitem{Lo} L.\ Lov\'asz, {\it Kneser's conjecture, chromatic number,
and homotopy}, J.\ Combin.\ Theory Ser.~A {\bf 25}, (1978), no.~3,
pp.~319--324.

\bibitem{Qu} D.~Quillen, {\em Higher algebraic K-theory} I, Lecture
Notes in Mathematics {\bf 341}, (1973), pp.~77--139, Springer-Verlag.

\bibitem{Zi02} G.M.~Ziegler, {\it Generalized Kneser coloring theorems
    with combinatorial proofs}, Invent.~Math. {\bf 147} (2002),
  pp.~671-691.


\end{thebibliography}
\end{document}